\documentstyle[11pt,amssymb,amstex]{amsart}

\newtheorem{thm}{Theorem}[section]
\newtheorem{lem}[thm]{Lemma}

\def\a{\alpha}
\def\o{\omega}
\def\O{\Omega}

\begin{document}

\title[Vector valued martingale-ergodic processes]
{Vector valued martingale-ergodic and ergodic-martingale theorems}

\maketitle
\begin{center}
\author{ F.A. Shahidi $^{\rm a}$ \footnote{Corresponding author. Email: farruh.shahidi@@gmail.com} and I.G. Ganiev $^{\rm b}$ }\\\vspace{12pt}
\address{$^{\rm a}${\em{Faculty of Information and Communication Technology,
International Islamic University Malaysia, P.O Box 10, 50728,
Kuala Lumpur, Malaysia}};

$^{\rm b}${\em{Faculty of Engineering, International Islamic
University Malaysia, P.O Box 10, 50728, Kuala Lumpur, Malaysia}}}

\end{center}

\begin{abstract} We prove martingale-ergodic and
ergodic-martingale theorems for vector valued Bochner integrable
functions. We obtain dominant and maximal inequalities. We also
prove weighted and multiparameter  martingale-ergodic and ergodic
martingale theorems.

\end{abstract}

\section{Introduction}
Since S. Kakutani \cite{Kak}, it has been known the similarity of
ergodic averages and martingales in terms of their behavior. In
addition, it was noticed that the proofs of ergodic theorems and
(reversed) martingale convergence theorems have some similarities.
Therefore the problem of unification of these processes was
naturally arisen. Since then several approaches have been
suggested by M. Jerrison (1959), G.-C. Rota (1961), A. and C.
Ionescu-Tulcea (1963), A.M. Vershik (1960s) and A.G. Kachurovskii
(see \cite{Kach2} for review). For example, M Jerrison showed that
ergodic averages can be considered as  martingales in some space
with $\sigma-$ finite measure. G.-C. Rota introduced so called
generalized martingales which allow to unify martingales and
ergodic averages with respect to Abel summation and proved the
convergence theorems for these processes. Abstract theorem of
Ionescu-Tulcea gives a unique proof of martingale and ergodic
convergence theorems. A.M. Vershik's approach was based on
consideration of actions of locally finite groups. Unfortunately,
each of the above approaches has some weaknesses by means of
unification of these two processes into a unique superstructure
\cite{Kach2}. Recently, A.G. Kachurovskii in \cite{Kach1}, and in
more detail in \cite{Kach2} developed a new theory of
martingale-ergodic and ergodic-martingale processes and proved
convergence theorems for these processes from which Birkhoff's
pointwise ergodic theorem, von Neumann's mean ergodic theorem
\cite{kren} as well as  Doob's martingale convergence theorems
\cite{Dub1}, can be obtained as the degenerate cases. Moreover,
the continuous analogues of martingale-ergodic and
ergodic-martingale processes were studied by I.V. Podvigin
\cite{pod1}, \cite{pod2}.

Note the  vector valued ergodic theory is developed sufficiently
well \cite{kren}, \cite{Chac}. Moreover, vector valued analogues
of Doob's martingales convergence theorems are also known and can
be found in \cite{Diestel}, \cite{Vaxan}. So, it is natural to
develop Kachurovskii's unification theory for vector valued
martingale and ergodic processes.

 The aim of this paper is to extend martingale-ergodic and
ergodic martingale processes to the space $L_p(\O, X)$ of Bochner
integrable functions with values in a Banach space $X.$ We prove
norm convergence as well as a.e. convergence for $X-$ valued
martingale-ergodic and ergodic-martingale processes. We establish
dominant and maximal inequalities. Also, we extent obtained
results for weighted and multiparameter martingale-ergodic and
ergodic martingale processes.

The paper is organized in the following way. In the next section
we define the vector valued martingale-ergodic processes and prove
convergence theorems for them. In the following section we
consider ergodic-martingale processes and prove convergence
theorems for them. Finally, the last section is devoted to
weighted and multiparameter extensions of obtained results.

\section{Vector valued martingale ergodic theorems}
In this section we prove a vector valued martingale ergodic
theorems.

 Throughout this paper $X$ will be denoted a reflexive Banach space
with the norm $||\cdot||_X$ and $(\O ,\beta, \mu)$ a finite
measure space. By $L_p(X)=L_p(\O, X), \ 1\le p<\infty$ we denote
the Banach space of X valued measurable functions $f$ on $\O.$
with the norm defined as

$$||f||_p=\left(\int_{\O}||f(\o)||_X^pd\mu\right)^{\frac 1p}.$$

We just write $L_p$ when $X=R.$ By $E(f|F)$ we denote the
conditional expectation of $f\in L_p(\O, X)$ with respect to
$\sigma-$ subalgebra $F$ of $\beta.$  Let $F_n$ be a sequence of
monotonically increasing (decreasing) $\sigma-$ subalgebras such
that $F_n\uparrow F_{\infty}$ ($F_n\downarrow F_{\infty}$) as
$n\to\infty.$ Henceforth we only consider a monotone sequence of
$\sigma-$ subalgebras $F_n.$ The sequence $(f_n)_{n\ge 1}$ in
$L_p(\O, X), \ 1\le p<\infty$ is called an \textit{ordinary
martingale}, if $f_n=E(f_{n+1}|F_n),$ and it is called a
\textit{reversed martingale}, if $f_{n+1}=E(f_n|F_n).$ A
martingale is called \textit{regular} if
$$f_n=E(f_0|F_n).$$
Henceforth we consider only regular martingales.

 We start with the following lemmas.

\begin{lem} Let $f_n$ be a sequence in $L_p(\O, X),\ p\ge 1.$
Assume that $f_n$ is norm convergent to $f^*$ in $L_p(\O, X)$ as
$n\to\infty.$ Then $E(f_{n_1}|F_{n_2})\rightarrow
E(f^*|F_{\infty})$ in norm on $L_p(\O, X)$ as $n_1,\
n_2\to\infty$(independently).
\end{lem}
\begin{pf}
 The proof can be given analogously as Lemma 1 of \cite{Kach2}.

\end{pf}

We say that a linear operator $T$ in $L_1(\O, X)$ is
\textit{positively dominated} if there exists a positive linear
contraction $T'$ in $L_1,$ called a \textit{positive dominant} of
$T,$ such that
$$||Tf||_X\le T'(||f||_X).$$

Let us now consider some examples of positively dominated
operators. More details can be found in \cite{kus2}, \cite{SuchF}.

1. If $T$ is a real valued, then it is positively dominated by
some positive linear contraction on $L_1.$ For the vector valued
$T$, a positive dominant may not exist in general.

2. Let $\tau$ be an endomorphism on $(\O ,\beta, \mu)$ then the
linear operator $T:L_1(\O, X)\rightarrow L_1(\O, X)$ given by
$Tf=f\circ\tau$ is said to be generated by $\tau.$ $T$ is
positively dominated by $T'$ with $T'(||f||_X)=||f||_X\circ\tau.$

3. Assume that the Banach space $X$ has Radon-Nikodym property
(this holds surely if $X$ is reflexive). Consider the conditional
expectation $E(f|F)$ with respect to $\sigma-$ subalgebra $F$ of
$\beta.$ For $f\in L_1(\O, X),$ the conditional expectation
$E(f|F)$ is Radon-Nikodym density with respect to the finite
measure $\mu$ on $F.$ Since $||E(f|F)||_X\le E'(||f||_X|F)$ a.e.
for all $f\in L_1(\O, X),$ where $E'(\cdot|F)$ is a conditional
expectation on $L_1,$ then the operator $E(\cdot|F)$ is positively
dominated by $E'(\cdot|F).$

\begin{lem}
Let $f_n$ be a sequence in $L_1(\O, X)$ and $f_n\rightarrow f^*$
a.e. as $n\to\infty.$ Assume that $h(\o)=sup_n||f_n(\o)||_X\in
L_1.$  Then $E(f_{n_1}|F_{n_2})\rightarrow E(f^*|F_{\infty})$ a.e.
as $n_1,\ n_2\to\infty$(independently).

\end{lem}

\begin{pf} Note that
$$||E(f_{n_1}|F_{n_2})-E(f^*|F_{\infty})||_X\le ||E(f_{n_1}|F_{n_2})-E(f^*|F_{n_2})||_X+||E(f^*|F_{n_2})-E(f^*|F_{\infty})||_X.$$
According to martingale convergence theorem\cite{Chatt} (see also
\cite{Pis1}) $||E(f^*|F_{n_2})-E(f^*|F_{\infty})||_X$ converges
a.e. to $0$  as $n_2\to\infty.$ Let us estimate
$||E(f_{n_1}|F_{n_2})-E(f^*|F_{n_2})||_X.$ We have
$$||E(f_{n_1}|F_{n_2})-E(f^*|F_{n_2})||_X=||E(f_{n_1}-f^*|F_{n_2})||_X\le$$
$$\le E'(||f_{n_1}-f^*||_X\ |F_{n_2})\le E'(h_{n_1}|F_{n_2})$$
where $h_{n_1}=sup_{m\ge n_1}||f_m-f^*||_X$ and $E'$ is the
positive dominant of $E$, that is, a conditional expectation on
$L_1.$. Since $f_n\rightarrow f^*$ a.e., then $h_{n_1}\to 0$ and
$h_{n_1}\le 2h\in L_1.$ Applying Lemma 2 of \cite{Kach1} we get
$E'(h_{n_1}|F_{n_2})\rightarrow 0.$ Therefore,
$E(f_{n_1}|F_{n_2})\rightarrow E(f^*|F_{\infty})$ a.e. as $n_1,\
n_2\to\infty$(independently).

\end{pf}

Let us put $$S_nf=\frac1n\sum\limits_{n=0}^{n-1}T^if,\
f^*=\lim\limits_{n\to\infty}S_nf,$$
$$S'_n(||f||_X)=\frac1n\sum\limits_{i=0}^{n-1}(T')^i(||f||_X).$$

\begin{thm}  Let $T$ be a linear operator generated by an endomorphism in $(\O ,\beta, \mu).$

1. Then
$$E(S_{n_1}f|F_{n_2})\rightarrow E(f^*|F_{\infty})$$
in $L_p(\O, X)$ as $n_1,\ n_2\to\infty,$ if $f\in L_p(\O, X),\
p\ge 1.$ Moreover $||E(f^*|F_{\infty}||_p\le ||f||_p.$

2. Let $sup_n||S_nf||_X\in L_1$  and $f\in L_p(\O, X).$ Then
$$E(S_{n_1}f|F_{n_2})\rightarrow E(f^*|F_{\infty})$$ a.e. in $L_1(\O,
X)$ as $n_1,\ n_2\to\infty,$ moreover $E(E(f^*|F_{\infty}))=Ef.$

\end{thm}

\begin{pf} The norm convergence comes from Lemma 2.1 and vector
valued analogue of mean ergodic theorem. Since the conditional
expectation operator is contracting in $L_p(\O, X),$ then
$||E(f^*|F_{\infty})||_p\le ||f^*||_p.$ Further, since $T$ is
generated by an endomorphism, then the ergodic average $S_nf$ is
contracting in $L_p(\O, X)$ norm. Hence $||f^*||_p\le ||f||_p.$

The convergence a.e. follows vector valued analogue of Birkhoff's
theorem (Theorem 4.2.1 in \cite{kren}) and Lemma 2.2.

\end{pf}

\textbf{Remark 1.} In the degenerate case when $F_{n_2}\equiv F,$
the above theorem coincides with the mean and a.e. vector valued
ergodic theorems. Moreover, if $T\equiv id$ the theorem coincides
with a vector valued martingale convergence theorem of S.
Chatterjee (see \cite{Chatt}, \cite{Diestel}, \cite{Vaxan}).

\textbf{Remark 2.} In case when $F_n$ be a sequence of
monotonically decreasing $\sigma-$ subalgebras such that
$F_n\downarrow F$ as $n\to\infty,$ Theorem 2.3 remains true.

\textbf{ Remark 3.} The condition "$sup_n||S_nf||_X\in L_1$" is
crucial. According to \cite{ArgRos}, this condition can not be
omitted even in real valued case. However, in real valued case,
under the assumption that conditional expectation and ergodic
average commute, the convergence theorem is given for the unified
processes in \cite{pod2}  without integrability of supremum.

Let $L(\O, X)[log^+L(\O, X)]^m $ denote the class of all functions
for which
 $$\int\limits_{X}||f||_X[log\ max(1, ||f||_X)]^m<\infty.$$
This space is a Banach space if we define the norm of $f$ as above
integral. We just write $L[log^+L]^m$ when $X=R.$

\begin{thm}(Dominant inequality) Let $T$ be a linear operator generated by
an endomorphism in $(\O ,\beta, \mu).$ Let $F_{n_2}\downarrow
F_{\infty}$ and $sup_n||S_nf||_X\in L_1.$  Then for $p>1,$ we have

1. $$||sup_{n_1,n_2}||E(S_{n_1}f|F_{n_2})||_X||_p\le \left(\frac
p{p-1}\right)^2||f||_p.$$

2. If $f\in L(\O, X)[log^+L(\O, X)]^{m+2},$ then
$$sup_{n_1,n_2}||E(S_{n_1}f|F_{n_2})||_X\in L[log^+L)]^m .$$
\end{thm}

\begin{pf} Let $g=sup_{n_1}||S_{n_1}f||_X.$ Then
$$||sup_{n_1,n_2}||E(S_{n_1}f|F_{n_2})||_X||_p\le
||sup_{n_2}E'(g|F_{n_2})||_p$$ where $E'$ is a positive dominant
of $E.$ Since $g$ is a real valued function, then by dominant
inequality for the reversed martingales \cite{Sh1}, we get
$$||sup_{n_2}E'(g|F_{n_2})||_p\le \frac p{p-1}||E'(g|F_1)||_p.$$

Since the conditional expectation is contracting, then
$$\frac p{p-1}||E'(g|F_1)||_p\le \frac p{p-1}||g||_p\le $$
$$\le \frac p{p-1}||sup_n||S_nf||_X||_p\le\frac p{p-1}||sup_nS'_n(||f||_X)||_p.$$

Now, again by using the dominant inequality \cite{kren}
$$||sup_nS'_n(||f||_X)||_p\le  \frac p{p-1}||||f||_X||_p= \frac p{p-1}||f||_p.$$
Therefore
$$||sup_{n_1,n_2}||E(S_{n_1}f|F_{n_2})||_X||_p\le \left(\frac
p{p-1}\right)^2||f||_p.$$

2. Let $g=sup_{n_1}||S_{n_1}f||_X.$ First we apply Theorem 2 from
\cite{Y}, which asserts that if the function $f$ is from the class
$ L(\O, X)[log^+L(\O, X)]^{m+2},$ then $g=sup_{n_1}||S_{n_1}f||_X$
belongs to the class $ L[log^+L]^{m+1}.$ We have the following
estimation

$$sup_{n_1,n_2}||E(S_{n_1}f|F_{n_2})||_X\le
sup_{n_2}E'(g|F_{n_2})$$ where $E'$ is a positive dominant of $E.$

Note that the sequence $E'(g|F_{n_2})_{n_2}$ is just the reversed
nonnegative submartingale. Since $g\in L[log^+L]^{m+1},$ then by
applying corollary 1.4 from \cite{Such}, we get
$sup_{n_2}E'(g|F_{n_2})\in L[log^+L]^m .$ So the above inequality
implies
$$sup_{n_1,n_2}||E(S_{n_1}f|F_{n_2})||_X\in L[log^+L)]^m .$$

\end{pf}

\begin{thm}(Maximal Inequality) Let $T$ be a linear operator generated by an endomorphism in $(\O ,\beta, \mu)$
and $f\in L_p(\O, X),\ p>1.$ Then the following inequality holds
for any $\varepsilon>0$ in the case when $F_n\downarrow
F_{\infty}$ and $sup_n||S_nf||_X\in L_1.$
$$\mu\{sup_{n_1,n_2}||E(S_{n_1}f|F_{n_2})||_X\ge\varepsilon\}\le \left(\frac p{p-1}\right)^p\frac{||f||_p^p}{\varepsilon^p}.$$

\end{thm}

\begin{pf} We put $g=sup_{n_1}||S_{n_1}f||_X.$ Then

$$\mu\{sup_{n_1,n_2}||E(S_{n_1}f|F_{n_2})||_X\ge\varepsilon\}\le\mu\{sup_{n_2}E'(g|F_{n_2})\ge\varepsilon\}$$
where $E'$ is a positive dominant of $E.$

By maximal inequality for reversed martingales \cite{Sh1},  we
have

$$\mu\{sup_{n_2}E'(g|F_{n_2})\ge\varepsilon\}\le \frac 1{\varepsilon^p}||E'(g|F_1)||_p^p.$$
Now, we apply the fact that the conditional expectation is a
contraction in $L_p$ and ergodic average $S_{n_1}$ is positively
dominated, we get
$$\frac 1{\varepsilon^p}||E'(g|F_1)||_p^p\le\frac 1{\varepsilon^p}||g||_p^p\le \frac 1{\varepsilon^p}||sup_{n_1}S'_{n_1}(||f||_X)||_p^p.$$
Further from dominant inequality for ergodic averages, we obtain

$$\frac 1{\varepsilon^p}||sup_{n_1}S'_{n_1}(||f||_X)||_p^p\le \left(\frac p{p-1}\right)^p\frac{||f||_p^p}{\varepsilon^p}.$$
That is why

$$\mu\{sup_{n_1,n_2}||E(S_{n_1}f|F_{n_2})||_X\ge\varepsilon\}\le \left(\frac p{p-1}\right)^p\frac{||f||_p^p}{\varepsilon^p}.$$

\end{pf}

We say that $T$ is an $L_1-L_{\infty}$ contraction if $||Tf||_1\le
||f||$ and $||Tf||_{\infty}\le ||f||_{\infty}$, where
$$||f||_1=\int\limits_{\O}||f(\o)||_X d\mu$$ and
$$||f||_{\infty}=inf\{\lambda: ||f(\o)||_X\le\lambda \ a.e\}.$$

Chacon \cite{Chac} proved the individual ergodic theorem for
$L_1-L_{\infty}$ contractions acting in $L_1(\O, X),$ where $X$ is
reflexive. The following theorem, which can be proved analogously
as Theorem 2.3, unifies Chacon's theorem and martingale ergodic
theorem.

\begin{thm} Let $T:L_1(\O,X)\rightarrow L_1(\O,X)$ be  positively dominated by an $L_1-L_{\infty}$ contraction $T'$ in $L_1.$
Then

1. $E(S_{n_1}f|F_{n_2})$ converges in $L_p(\O,X),\ p\ge 1,$ if
$f\in L_p(\O, X)$ as $n_1,n_2\to\infty.$

2. $E(S_{n_1}f|F_{n_2})$ converges a.e. if
$sup_{n_1}||S_{n_1}f||_X\in L_1.$
\end{thm}

\begin{pf} Note that since $T'$ is an $L_1-L_{\infty}$ contraction in
$L_1,$ then $T$ is also an $L_1-L_{\infty}$ contraction in
$L_1(\O,X).$ Therefore,  both assertions are the consequences of
Chacon's theorem \cite{Chac}, Lemmas 2.1 and 2.2.

\end{pf}

The following dominant and maximal inequalities for
$L_1-L_{\infty}$ contraction hold with the same constants of
Theorems 2.4 and 2.5.

\begin{thm} Under the hypothesis of Theorem 2.6 the following
inequalities hold true if $sup_{n_1}||S_{n_1}f||_X\in L_1$ and
$p>1.$

1.
$$||sup_{n_1,n_2}||E(S_{n_1}f|F_{n_2})||_X||_p\le \left(\frac
p{p-1}\right)^2||f||_p.$$

2.$$\mu\{sup_{n_1,n_2}||E(S_{n_1}f|F_{n_2})||_X\ge\varepsilon\}\le
\left(\frac p{p-1}\right)^p\frac{||f||_p^p}{\varepsilon^p}.$$

\end{thm}

\section{Vector valued Ergodic-Martingale processes}
In this section, we define vector valued ergodic-martingale
processes and prove the convergence theorems for them.

\begin{lem} Let $T$ a linear operator generated by an endomorphism $(\O,\beta, \mu)$ and
$f_n$ be a sequence in $L_p(\O, X),\ p\ge
1,$ and $f_n\rightarrow f$ in $L_p(\O, X)$ whenever $n\to\infty.$
Then
$$S_{n_1}f_{n_2}\rightarrow
f^*=\lim\limits_{n\to\infty}S_nf$$ in $L_p$ as $n_1,n_2\to\infty.$
\end{lem}

\begin{pf} Since
$$||S_{n_1}f_{n_2}-f^*||_p\le ||S_{n_1}(f_{n_2}-f)||_p+||S_{n_1}f-f^*||_p,$$
then from condition of the lemma and by virtue of contraction
property of ergodic averaging $S_n$ we get the desired result.
\end{pf}

The following lemma is a generalization of Theorem 7.5 proved by
Maker P.T.(see \cite{kren}, p.66 or \cite{Kach2}), in a vector
valued case.

\begin{lem} Let $T$ a linear operator generated by an endomorphism $(\O,\beta,
\mu)$ and $f_n$ be a sequence  in $L_p(\O, X),\ p\ge 1$ with
$f_n\rightarrow f$ a.e. as $n\to\infty$ in $||\cdot||_X$ norm,
$h=sup_n||f_n||_X\in L_1.$  Then for positively dominated operator
$T$ we have
$$S_{n_1}f_{n_2}\rightarrow
f^*=\lim\limits_{n\to\infty}S_nf$$ a.e. as $n_1,n_2\to\infty$
\end{lem}

\begin{pf}One can see that

$$||S_{n_1}f_{n_2}-f^*||_X\le ||S_{n_1}(f_{n_2}-f)||_X+||S_{n_1}f-f^*||_X$$
Obviously, $||S_{n_1}(f)-f^*||_X\rightarrow 0$ a.e.

Since $T$ is positively dominated, then there exists a positive
dominant $T'$ such that
$$||Tf||_X\le T'(||f||_X).$$
Hence
$$||S_{n_1}(f_{n_2}-f)||_X\le
||\frac1{n_1}\sum\limits_{i=0}^{n_1-1}T^i(f_{n_2}-f)||_X\le$$
$$\le\frac1{n_1}\sum\limits_{n=0}^{n_1-1}(T')^i(||f_{n_2}-f||_X)=S'_{n_1}(||f_{n_2}-f||_X).$$

Note that $||f_{n_2}-f||_X\in L_1$ and according to our assumption
$||f_{n_2}-f||_X\rightarrow 0$ a.e. Further, since the function
$h=sup_n||f_n||_X$ is integrable, then due to $||f_{n_2}-f||_X\le
2h$ the function $||f_{n_2}-f||_X$ is also integrable. According
to Theorem 7.5 \cite{kren} (p.66) we get
$S'_{n_1}(||f_{n_2}-f||_X)\rightarrow 0$ a.e. Therefore,
$||S_{n_1}f_{n_2}-f^*||_X\rightarrow 0$ a.e.

\end{pf}

\begin{thm} Let $T$ be a linear operator generated by an endomorphism $(\O,\beta, \mu)$.
The following statements hold.

1.
$$S_{n_1}(E(f|F_{n_2}))\rightarrow E(f^*|F_{\infty})$$
in $L_p(\O, X)$ as $n_1,\ n_2\to\infty,$ if $f\in L_p(\O, X),\
p\ge 1.$ Moreover $||E(f^*|F_{\infty})||_p\le ||f||_p.$

2. If $sup_{n_2}||E(f|F_{n_2})||_X\in L_1,$  then
$$S_{n_1}(E(f|F_{n_2}))\rightarrow f_{\infty}^*$$ a.e. in $L_1(\O,
X)$ as $n_1,\ n_2\to\infty,$ moreover $E(E(f^*|F_{\infty}))=Ef.$

\end{thm}

\begin{pf} The proofs are the same as the proofs of Theorem 2.3.
The norm convergence is the combination of Lemma 3.1, vector
valued mean ergodic theorem \cite{kren}, and vector valued
martingale convergence theorem. Note that if a linear operator is
generated by an endomorphism, then it possesses a positive
dominant. So the a.e. convergence is the combination of Lemma 3.2,
vector valued individual ergodic theorem and vector valued
martingale convergence theorem.
\end{pf}

 As the matter of fact, as in classical $L_p$ spaces, dominant
and maximal inequalities for ergodic-martingale processes can be
obtained by the same arguments of Theorems 2.4 and 2.5 for
martingale-ergodic processes. However, as in real valued case
\cite{Kach2}, the advantage of ergodic-martingale processes is
that the sequence of martingales need not necessarily to be
reversed.

\begin{thm}(Dominant Inequality) For a linear operator $T$ generated by an endomorphism  on $(\O,\beta, \mu)$
the following estimates hold true.

1. If $f\in L_p(\O, X),\ p>1,$ and $sup_{n_2}||E(f|F_{n_2})||_X\in
L_1,$ then
$$||sup_{n_1,n_2}||S_{n_1}(E(f|F_{n_2}))||_X||_p\le \left(\frac
p{p-1}\right)^2||f||_p.$$

2. If $f\in L(\O, X)[log^+L(\O, X)]^{m+2},$ then
$$sup_{n_1,n_2}||S_{n_1}(E(f|F_{n_2}))||_X\in L[log^+L]^m .$$
\end{thm}

\begin{pf} The proof of this Theorem is very similar to the proof
of Theorem 2.4. Let $g=sup_n||E(f|F_n)||_X.$ Then

$$||sup_{n_1,n_2}||S_{n_1}(E(f|F_{n_2}))||_X||_p\le||sup_{n_1}S'_{n_1}g||_p.$$

From the dominant inequality for ergodic averages \cite{kren} we
get

$$||sup_{n_1}S'_{n_1}g||_p\le \frac p{p-1}||g||_p.$$

Since the conditional expectation operator is positively
dominated, then

$$\frac p{p-1}||g||_p=\frac p{p-1}||sup_n||E(f|F_n)||_X||_p\le
\le \frac p{p-1}||sup_nE(||f||_X|F_n)||_p\left(\frac
p{p-1}\right)^2||f||_p.$$

Now we prove part 2. We set $g=sup_n||E(f|F_n)||_X.$ One has
$$g\le sup_{n_2}E'(||f||_X|F_{n_2}).$$

Note that since $f\in L(\O, X)[log^+L(\O, X)]^{m+2},$ then
$||f||_X\in L[log^+L]^{m+2}.$ Since the sequence
$E'(||f||_X|F_{n_2})$ is reversed nonnegative martingale, then
Corrolary 1.4 of \cite{Such} implies
$$g=sup_n||E(f|F_n)||_X\in L[log^+L]^{m+1}.$$
Further, note that
$$sup_{n_1,n_2}||S_{n_1}(E(f|F_{n_2}))||_X\le sup_{n_1}S'_{n_1}g.$$

 Moreover, the condition $g\in L[log^+L]^{m+1}$ yields
$S'_ng\in L[log^+L]^m$ \cite{kren}. Finally, the above inequality
implies

$$sup_{n_1,n_2}||S_{n_1}(E(f|F_{n_2}))||_X\in L[log^+L]^m .$$

\end{pf}

\begin{thm}(Maximal Inequality) Let $T$ be a linear operator generated by an endomorphism $(\O,\beta, \mu)$
If $sup_{n_2}||E(f|F_{n_2})||_X\in L_1$, then for any
$\varepsilon>0$ we have

$$\mu\{sup_{n_1,n_2}||S_{n_1}(E(f|F_{n_2}))||_X\ge\varepsilon\}\le \left(\frac p{p-1}\right)^p\frac{||f||_p^p}{\varepsilon^p}.$$

\end{thm}

\begin{pf} $g=sup_n||E(f|F_n)||_X.$ Then

$$\mu\{sup_{n_1,n_2}||S_{n_1}(E(f|F_{n_2}))||_X\ge\varepsilon\}\le
\mu\{sup_{n_1}S'_{n_1}g\ge\varepsilon\}.$$

By the maximal inequality for ergodic averages \cite{kren}, we
obtain

$$\mu\{sup_{n_1}S'_{n_1}g\ge\varepsilon\}\le\frac 1{\varepsilon^p}||g||_p^p.$$

Finally, from the dominant inequality for vector valued
martingales \cite{Pis1}, \cite{Pis2}, we get

$$\frac 1{\varepsilon^p}||g||_p^p=\frac 1{\varepsilon^p}||sup_n||E(f|F_n)||_X||_p^p\le \left(\frac p{p-1}\right)^p\frac{||f||_p^p}{\varepsilon^p}, $$
hence

$$\mu\{sup_{n_1,n_2}||S_{n_1}(E(f|F_{n_2}))||_X\ge\varepsilon\}\le \left(\frac p{p-1}\right)^p\frac{||f||_p^p}{\varepsilon^p}.$$

\end{pf}

\begin{thm} Let an operator $T:L_1(\O,X)\rightarrow L_1(\O,X)$ be positively dominated by an
$L_1-L_{\infty}$ contraction $T'$ in $L_1.$ Then

1. $S_{n_1}(E(f|F_{n_2}))$ converges in $L_p(\O,X),\ p\ge 1,$ if
$f\in L_p(\O, X)$ as $n_1,n_2\to\infty.$

2. $S_{n_1}(E(f|F_{n_2}))$ converges a.e. if
$sup_{n_2}||E_{n_2}(f|F_{n_2})||_X\in L_1$  and $f\in L_1(\O, X)$.
\end{thm}

\begin{pf} Since $T'$ is an $L_1-L_{\infty}$ contraction in
$L_1,$ then $T$ is also an $L_1-L_{\infty}$ contraction in
$L_1(\O,X).$ Therefore, the proof comes from Chacon's theorem
\cite{Chac}, Lemmas 3.1 and 3.2.

\end{pf}

\begin{thm}  Let an operator $T:L_1(\O,X)\rightarrow L_1(\O,X)$ be positively dominated by an
$L_1-L_{\infty}$ contraction in $L_1.$  If $f\in L_p(\O, X),\ p>1$
and $sup_{n_2}||E(f|F_{n_2})||_X\in L_1,$ then

1. $$||sup_{n_1,n_2}||S_{n_1}E(f|F_{n_2})||_X||_p\le \left(\frac
p{p-1}\right)^2||f||_p.$$

2.
$$\mu\{sup_{n_1,n_2}||S_{n_1}E(f|F_{n_2})||_X\ge\varepsilon\}\le
\left(\frac p{p-1}\right)^p\frac{||f||_p^p}{\varepsilon^p}.$$

\end{thm}

\section{Weighted and multiparameter cases}

In fact, not only can martingales be unified by ergodic averages,
but also multiparameter martingales can be unified by weighted and
multiparameter ergodic averages. In this section we provide a
weighted and multiparameter martingale-ergodic as well as
ergodic-martingale theorems.

Firstly, we define the following terminology. Let $\a_i,\ i\in N$
be a sequence of complex numbers. We say that $\a_i$ is a
\textit{bounded Besicovitch sequence} if the sequence $\a_i$ is
bounded and for any $\varepsilon>0$ there exists a trigonometric
polynomial $\varphi_{\varepsilon}$ such that
$$\lim\limits_{n\to\infty}\frac1n\sum\limits_{i=0}^{n-1}|\a_i-\varphi_{\varepsilon}(i)|<\varepsilon.$$

We define the following weighted average

$$S_n(T,\ \a,\ f)=\frac1n\sum\limits_{i=0}^{n-1}\a_iT^if,$$
where $T$ is a linear operator in $L_1(\O, X).$ When $X$ is
reflexive and $T$ is an $L_1-L_{\infty}$ contraction, then norm
and a.e. convergence of $S_n(T,\ \a,\ f)$ is due to K. Berdan
\cite{berdan}. The following theorem is a unification of vector
valued weighted ergodic theorem and martingale convergence
theorem.

\begin{thm} Let an operator $T:L_1(\O,X)\rightarrow L_1(\O,X)$ be positively dominated by an
$L_1-L_{\infty}$ contraction $T'$ in $L_1.$ Let $f\in L_1(\O, X)$
such that $sup_n||S_n(T,\ \a,\ f||_X\in L_1.$   Then

1. $E(S_{n_1}(T,\ \a,\ f)|F_{n_2})$ converges in $L_p(\O,X),\ p\ge
1,$ if $f\in L_p(\O, X)$ as $n_1,n_2\to\infty.$

2. $E(S_{n_1}(T,\ \a,\ f)|F_{n_2})$ converges a.e.

In case $F_{n_2}\downarrow F_{\infty}$ and $p>1$ we have

3.$$||sup_{n_1,n_2}||E(S_{n_1}(T,\ \a,\ f)|F_{n_2})||_X||_p\le
\a\left(\frac p{p-1}\right)^2||f||_p.$$

4.
$$\mu\{sup_{n_1,n_2}||E(S_{n_1}(T,\ \a,\ f)|F_{n_2})||_X\ge\varepsilon\}\le
\a\left(\frac p{p-1}\right)^p\frac{||f||_p^p}{\varepsilon^p},$$
where $\a=sup_i(\a_i).$

\end{thm}

\begin{pf} The proof of the first part is the combination of Lemma
2.1 and norm convergence of $S_n(T,\ \a,\ f)$ \cite{berdan}. The
second part comes from Lemma 2.2. and a.e. convergence of $S_n(T,\
\a,\ f)$ \cite{berdan}.

3. Let $g=sup_{n_1}||S_{n_1}(T,\ \a,\ f)||_X.$ Then
$$||sup_{n_1,n_2}||E(S_{n_1}(T,\ \a,\ f)|F_{n_2})||_X||_p\le
||sup_{n_2}E'(g|F_{n_2})||_p,$$ where $E'$ is the positive
dominant of $E.$

 Since $g$ is a real valued
function, then by dominant inequality for the reversed martingales
we get
$$||sup_{n_2}E'(g|F_{n_2})||_p\le \frac p{p-1}||E'(g|F_1)||_p.$$

Since the conditional expectation is contracting, then
$$\frac p{p-1}||E'(g|F_1)||_p\le \frac p{p-1}||g||_p\le $$
$$\le \frac p{p-1}||sup_n||S_n(T,\ \a,\ f)||_X||_p\le\frac p{p-1}||sup_nS'_n(T', \a \ ||f||_X)||_p$$
where $S'_n(T', \a \
||f||_X)=\frac1n\sum\limits_{i=0}^{n-1}\a_i(T')^i(||f||_X).$

Now, by using the dominant inequality \cite{berdan} for $S'_n(T,\
\a,\ f)$, in case $X=R$, we get
$$||sup_nS'_n(||f||_X)||_p\le \a\frac p{p-1}||||f||_X||_p=\a\frac p{p-1}||f||_p.$$
Hence, we get
$$||sup_{n_1,n_2}||E(S_{n_1}(T,\ \a,\ f)|F_{n_2})||_X||_p\le \a\left(\frac
p{p-1}\right)^2||f||_p.$$

4.  We put $g=sup_{n_1}||S_{n_1}(T,\ \a,\ f)||_X.$ Then

$$\mu\{sup_{n_1,n_2}||E(S_{n_1}(T,\ \a,\ f)|F_{n_2})||_X\ge\varepsilon\}\le\mu\{sup_{n_2}E'(g|F_{n_2})\ge\varepsilon\}$$
where $E'$ is a positive dominant of $E.$

By maximal inequality for reversed martingales \cite{Sh1},  we
have

$$\mu\{sup_{n_2}E'(g|F_{n_2})\ge\varepsilon\}\le \frac 1{\varepsilon^p}||E'(g|F_1)||_p^p.$$
Now, we apply again that the conditional expectation is a
contraction in $L_p$ and get
$$\frac 1{\varepsilon^p}||E'(g|F_1)||_p^p\le\frac 1{\varepsilon^p}||g||_p^p\le \frac 1{\varepsilon^p}||sup_{n_1}S'_{n_1}(T,\ \a,\ ||f||_X)||_p^p$$
where $S'_n(T', \a \
||f||_X)=\frac1n\sum\limits_{i=0}^{n-1}\a_i(T')^i(||f||_X).$

Further from dominant inequality weighted for ergodic averages
\cite{berdan}, we obtain

$$\frac 1{\varepsilon^p}||sup_{n_1}S'_{n_1}(T', \a \
||f||_X)||_p^p\le \a\left(\frac
p{p-1}\right)^p\frac{||f||_p^p}{\varepsilon^p}.$$ That is why

$$\mu\{sup_{n_1,n_2}||E(S_{n_1}(T,\ \a,\ f)|F_{n_2})||_X\ge\varepsilon\}\le \a\left(\frac p{p-1}\right)^p\frac{||f||_p^p}{\varepsilon^p}.$$

\end{pf}

The following theorem is a weighted vector valued
ergodic-martingale theorem.

\begin{thm}  Let an operator $T:L_1(\O,X)\rightarrow L_1(\O,X)$ be positively dominated by an
$L_1-L_{\infty}$ contraction $T'$ in $L_1.$ Let $f\in L_1(\O, X)$
such that $sup_n|| E(f|F_n)||_X\in L_1.$  Then

1. $S_{n_1}(T,\ \a,\ E(f|F_{n_2}))$ converges in $L_p(\O,X),\ p\ge
1,$ if $f\in L_p(\O, X)$ as $n_1,n_2\to\infty.$

2. $S_{n_1}(T,\ \a,\ E(f|F_{n_2}))$ converges a.e.

 In case when $F_{n_2}\downarrow F_{\infty}$ and $p>1$ we have

3.$$||sup_{n_1,n_2}||S_{n_1}(T,\ \a,\ E(f|F_{n_2}))||_X||_p\le
\a\left(\frac p{p-1}\right)^2||f||_p.$$

4. $$\mu\{sup_{n_1,n_2}||S_{n_1}(T,\ \a,\
E(f|F_{n_2}))||_X\ge\varepsilon\}\le \a\left(\frac
p{p-1}\right)^p\frac{||f||_p^p}{\varepsilon^p},$$ where
$\a=sup_i(\a_i).$

\end{thm}

\begin{pf} The proofs of 1 and 2 are the combinations of Lemma
3.1, 3.2 and weighted vector ergodic theorem of K.Berdan
\cite{berdan}. 3 and 4 can be obtained in a similar way to 3 and 4
of the previous theorem.
\end{pf}

Now, we turn to the unification of multiparmeter martingales and
ergodic averages. For each fixed $k,$ let $F_n^k, \ n\in N$ be
either increasing or decreasing $\sigma- $ subalgebras such that
$F_n^k\uparrow F_{\infty}^k$ or $F_n^k\downarrow F_{\infty}^k.$
Let $E_n^k=E(\cdot| F_n^k),\ n\in N$ and
$E_s=E_{s_1}^1E_{s_2}^2\cdots E_{s_{p+1}}^{p+1},$ for
$s=(s_1,s_2,\cdots, s_{p+1})\in N^{p+1}.$ By $(F_s,\ s\in
N^{p+1})$ we denote the increasing or decreasing net of $\sigma- $
subalgebras.

Consider the following multiparameter weighted average
$$S(T_d,\ \a_d,\  n_d,\ f)=\frac1{n_1n_2\cdots n_d}\sum\limits_{k_i=0 \ i=\overline{1,d}}^{n_i-1}\a_{k_1}^1\cdots
\a_{k_d}^dT_1^{k_1}\cdots T_d^{k_d}(f)$$ where $\{(\a_n^j)\}$ are
bounded Besicovitch sequences, $j=1,2,\cdots, d.$ We set
$\a=sup_{k=1,2,\cdots, d}sup_j|\a_k^j|.$

We define the multiparameter martingale-ergodic average as
$E_s(S(T_d,\ \a_d,\ n_d\ f)|F_s).$

\begin{thm}  Let the operators $T_i:L_1(\O,X)\rightarrow L_1(\O,X), \ i=\overline{1,d}$
be positively dominated by an $L_1-L_{\infty}$ contractions $T'_i$
in $L_1.$ Assume $f\in L_1(\O, X)$ and $sup_{n_j}||S_{n_j}(T,\
\a_d,\ n_d,\ f)||_X$ is integrable, then

1. The multiparameter martingale-ergodic average as $E_s(S(T_d,\
\a_d,\  n_d\ f)|F_s)$ converges a.e. as $n_d, s\to\infty$
independently.

For $f\in L_p(\O, X),\ p>1, $ we have

2.$$||sup_{n_j,s}||E_s(S(T_d,\ \a_d,\ n_d,\ f)|F_s)||_X||_p\le
\a\left(\frac p{p-1}\right)^{d+p+1}||f||_p.$$

3.$$\mu\{sup_{n_j,s}||E(S(T_d,\ \a_d,\ n_d,\
f)|F_s)||_X\ge\varepsilon\}\le \a^p\left(\frac
p{p-1}\right)^{pd}\frac{||f||_p^p}{\varepsilon^p}.$$

\end{thm}

\begin{pf} 1. Since the operators $T_i$ are dominated by an $L_1-L_{\infty}$ contractions
$T'_i,$ then $T_i$ are also  $L_1(\O, X)-L_{\infty}(\O, X)$
contractions. Therefore, $T_i$ satisfy the conditions of Theorem 4
of \cite{berdan2}, from which it follows that $S(T_d,\ \a_d,\
n_d,\ f)$ converges a.e. as $n_j\to\infty, \ j=1,2,\cdots d.$

On the other hand, since $X$ is reflexive, then it possesses the
Radon-Nykodim property, therefore, from Theorem 6.2 of
\cite{SuchF}, it follows that for $f\in L_p(\O, X)$
$$E_s=E_{s_1}^1E_{s_2}^2\cdots E_{s_{p+1}}^{p+1}(f)$$
converges a.e. as the indices $s_i\to\infty$ independently.

These two facts along with Lemma 2.2 prove the assertion 1 of the
theorem.

2. Let $g=sup_{n_j}||S(T_d,\ \a_d,\ n_d,\ f)||_X.$ Then we have

$$||sup_{n_j,s}||E_s(S(T_d,\ \a_d,\ n_d,\ f)|F_s)||_X||_p\le ||sup_{s}||E'_s(g|F_s)||_X||_p,$$
where $E'_s$ is a positive dominant of $E_s.$

Dominant inequality for reversed martingales and simple iteration
argument (see \cite{Pis1}, p.22) imply

$$||sup_{s}||E'_s(g|F_s)||_X||_p\le \left(\frac p{p-1}\right)^{p+1}||g||_p.$$

Further,
$$(\frac p{p-1})^{p+1}||g||_p=\left(\frac p{p-1}\right)^{p+1}||sup_{n_j}||S(T_d,\ \a_d,\ n_d,\
f)||_X||_p\le$$$$
 \le\left(\frac p{p-1}\right)^{p+1}||sup_{n_j}||S'(T'_d,\ \a_d,\ n_d,\ ||f||_X)||_X||_p,$$
where $S',$ given by
$$S'(T'_d,\ \a_d,\  n_d,\ ||f||_X)=\frac1{n_1n_2\cdots n_d}\sum\limits_{k_i=0 \ i=\overline{1,d}}^{n_i-1}\a_{k_1}^1\cdots
\a_{k_d}^d(T'_1)^{k_1}\cdots (T'_d)^{k_d}(g)$$ is a positive
dominant of $S.$

Finally, from  weighted multiparameter dominant inequality for
ergodic average $S'$ \cite{berdan2} (in case $X=R$), we have
$$||sup_{n_j}||S'(T'_d,\ \a_d,\ n_d,\ ||f||_X)||_p\le \a\left(\frac p{p-1}\right)^d||f||_p$$
and hence

$$||sup_{n_j,
s}||E_s(S(T_d,\ \a_d,\ n_d,\ f)|F_s)||_X||_p\le \a\left(\frac
p{p-1}\right)^{d+p+1}||f||_p.$$

3. Let $g=sup_{n_j}||S(T,\ \a_d,\ n_d,\ f)||_X$. Then
$$\mu\{sup_{n_j,s}||E(S(T_d,\ \a_d,\ n_d,\
f)|F_s)||_X\ge\varepsilon\}\le\mu\{sup_{s}E'_s(g|F_s)\ge\varepsilon\},$$
where each $(E^i_{s_i})'\ i=\overline{1,p+1}$ is the positive
dominant of $E^i_{s_i}\ i=\overline{1,p+1}.$

By maximal inequality for $(E^1_{s_1})'$ we have

$$\mu\{sup_{s}E'_s(g|F_s)\ge\varepsilon\}\le\frac 1{\varepsilon^p}
||(E^2_{s_2})'\cdots (E^{p+1}_{s_{p+1}})'(g|F_{s_2,\cdots,
s_{p+1}})||_p^p$$

From contraction property of conditional expectations, we have
$$\frac 1{\varepsilon^p}
||(E^2_{s_2})'\cdots (E^{p+1}_{s_{p+1}})'(g|F_{s_2,\cdots,
s_{p+1}})||_p^p\le \frac 1{\varepsilon^p}||g||_p^p.$$

Now we use the dominant inequality for ergodic averages and get
$$\frac 1{\varepsilon^p}||g||_p^p=\frac 1{\varepsilon^p}||sup_{n_j}||S(T,\ \a_d,\ n_d,\
f)||_X||_p^p\le$$
$$\le \frac 1{\varepsilon^p}\a^p\left(\frac p{p-1}\right)^{pd}||f||_p^p=
\a^p\left(\frac p{p-1}\right)^{pd}\frac
{||f||_p^p}{\varepsilon^p}$$

\end{pf}

\textbf{Remark 4.} It should be stressed that in above theorems
the number of conditional expectations depend on $p$ (there are
$p+1$ conditional expectations). That is why Theorem 4.1 can not
be considered as a particular case of Theorem 4.3. However, one
can unify one parameter martingale with weighted multiparamater
ergodic averages and obtain similar results.

The following is a multiparameter ergodic-martingale theorem.
However, maximal inequality for this process is unknown for us.

\begin{thm}  Let the operators $T_i:L_1(\O,X)\rightarrow L_1(\O,X),
\ i=\overline{1,d}$ be positively dominated by an $L_1-L_{\infty}$
contractions $T'_i$ in $L_1.$ Assume $f\in L_1(\O, X)$ and
$sup_s||E_s(f|F_s)||_X$ is integrable, then

1. The multiparameter ergodic-martingale average as $S(T_d,\
\a_d,\  n_d\, E_s(f|F_s))$ converges a.e. as $n_d, s\to\infty$
independently.

2. For $f\in L_p(\O, X),\ p>1. $ we have

$$||sup_{n_j,s}||S(T_d,\ \a_d,\ n_d,\ E_s(f|F_s))||_X||_p\le \a\left(\frac
p{p-1}\right)^{d+p+1}||f||_p.$$

\end{thm}

\begin{pf} 1. The proof is very similar to the proof of the previous
theorem.  Since the operators $T_i$ are dominated by an
$L_1-L_{\infty}$ contractions $T'_i,$ then $T_i$ are also
$L_1(\O, X)-L_{\infty}(\O, X)$ contractions. Therefore, $T_i$
satisfy the conditions of Theorem 4 of \cite{berdan2}, from which
it follows that $S(T_d,\ \a_d,\ n_d,\ f)$ converges a.e. as
$n_j\to\infty, \ j=1,2,\cdots d.$

Moreover, the reflexivity of $X$ allows us to use Theorem 6.2 of
\cite{SuchF}, which says that  for $f\in L_p(\O, X),$ the
multiparameter conditional expectations
$$E_s=E_{s_1}^1E_{s_2}^2\cdots E_{s_{p+1}}^{p+1}(f)$$
converges a.e. as the indices $s_i\to\infty$ independently.

Therefore, Lemma 3.2 concludes the proof of the assertion 1 of the
theorem.

2. Let  $g=sup_s||E_s(f|F_s)||_X.$ Then

$$||sup_{n_j,s}||S(T_d,\ \a_d,\ n_d,\ E_s(f|F_s))||_X||_p\le ||sup_{n_j}S'(T_d,\ \a_d,\ n_d,\ g )||_p$$
where $S'(T_d,\ \a_d,\ n_d,\ g),$ defined by
$$S'(T'_d,\ \a_d,\  n_d,\ f)=\frac1{n_1n_2\cdots n_d}\sum\limits_{k_i=0 \ i=\overline{1,d}}^{n_i-1}\a_{k_1}^1\cdots
\a_{k_d}^d(T'_1)^{k_1}\cdots (T'_d)^{k_d}(f)$$ is a positive
dominant of $S(T_d,\ \a_d,\ n_d,\ f).$

From dominant inequality for  weighted multiparameter ergodic
averages , we have

$$||sup_{n_j}S'(T_d,\ \a_d,\ n_d,\ g )||_p\le \a\left(\frac p{p-1}\right)^d||g||_p.$$

Note that since conditional expectations $E_s$ are positively
dominated, then

$$||E_s(f|F_s)||_X\le E'_s(||f||_X|F_s),$$ where
$E'_s$ is multiparameter real valued conditional expectation.
Thus,

$$\a\left(\frac p{p-1}\right)^d||g||_p=\a\left(\frac p{p-1}\right)^d||sup_s||E_s(f|F_s)||_X||_p\le \a\left(\frac p{p-1}\right)^d||E'_s(||f||_X|F_s)||_p.$$

The dominant inequality for reversed martingales and simple
iteration argument bring to the following estimate

$$\a\left(\frac p{p-1}\right)^d||E'_s(||f||_X|F_s)||_p\le \a\left(\frac p{p-1}\right)^{d+p+1}||f||_p.$$

That is why $$||sup_{n_j,s}||S(T_d,\ \a_d,\ n_d,\
E_s(f|F_s))||_X||_p\le \a\left(\frac
p{p-1}\right)^{d+p+1}||f||_p$$ holds.

\end{pf}

\end{document}